# Convex Matroid Optimization

Shmuel Onn *

**Abstract**

We consider a problem of optimizing convex functionals over matroid bases. It is richly expressive and captures certain quadratic assignment and clustering problems. While generally NP-hard, we show it is polynomial time solvable when a suitable parameter is restricted.

## 1 Introduction

Let $M = (N, \mathcal{B})$ be a matroid over $N := \{1, \ldots, n\}$ with collection of bases $\mathcal{B} \subseteq 2^N$. Let $w : N \longrightarrow \mathbb{R}^d$ be a weighting of matroid elements by vectors in $d$-space. For any subset $J \subseteq N$ let $w(J) := \sum_{j \in J} w(j)$ with $w(\phi) := 0$. Finally, let $c : \mathbb{R}^d \longrightarrow \mathbb{R}$ be a convex functional on $\mathbb{R}^d$. We consider the following algorithmic problem.

**Convex matroid optimization.** Given data as above, find a basis $B \in \mathcal{B}$ maximizing $c(w(B))$.

We begin with some examples of specializations of this problem.

**Example 1.1 Linear matroid optimization.** This is the special case of our problem with $d = 1$, $w : N \longrightarrow \mathbb{R}$ a weighting of elements by scalars, and $c : \mathbb{R} \longrightarrow \mathbb{R} : x \mapsto x$ the identity. The problem is to find a basis of maximum weight, and is quickly solvable by the greedy algorithm.

**Example 1.2 Positive semidefinite quadratic assignment.** This is the NP-hard problem [6] of finding a vector $x \in \{0, 1\}^n$ maximizing $||Wx||^2 = x^T W^T W x$ with $W$ a given $d \times n$ matrix. For fixed $d$ it is solvable in polynomial time [3]. The variant of this problem in which one asks for $x$ with restricted support $|\text{supp}(x)| = r$ is the special case of our problem with $M := U^r_n$ the uniform matroid of rank $r$ over $N$, with $w(j) := W^j$ the $j$th column of $W$ for all $j \in N$, and with $c : \mathbb{R}^d \longrightarrow \mathbb{R} : x \mapsto ||x||^2$ the $l_2$-norm (squared or not). The positive semidefinite quadratic assignment problem can be solved by solving the variant for $r = 0, \ldots, n$ and picking the best $x$.

---

*Research supported in part by a grant from ISF - the Israel Science Foundation, by a VPR grant at the Technion, and by the Fund for the Promotion of Research at the Technion.

**Example 1.3 Balanced clustering.** This is the problem of partitioning a given set $\{w_1, \ldots, w_n\}$ of points in $\mathbb{R}^d$ into two clusters $C_1, C_2$ of equal size $m := \frac{n}{2}$ so as to minimize the sum of cluster variances given by

$$\frac{1}{m} \sum_{w_j \in C_1} ||w_j - (\frac{1}{m} \sum_{w_j \in C_1} w_j)||^2 \quad + \quad \frac{1}{m} \sum_{w_j \in C_2} ||w_j - (\frac{1}{m} \sum_{w_j \in C_2} w_j)||^2 \ .$$

It can be shown by suitable manipulation of the variance expression that this is the special case of our problem with $M$ the uniform matroid of rank $\frac{n}{2}$ over $N$, with $w(j) := w_j$ for all $j \in N$, and with $c : \mathbb{R}^d \longrightarrow \mathbb{R} : x \mapsto ||x||^2 + ||w(N) - x||^2$ with $w(N) = \sum_{j=1}^{n} w_j$ the sum of all points.

While the linear matroid optimization problem (Example 1.1) is greedily solvable (cf. [4]), the general convex matroid optimization problem is NP-hard as indicated by Example 1.2. Nevertheless, in this article we show that, so long as $d$ is fixed, the problem can be solved in polynomial time for an arbitrary matroid $M$ and an arbitrary convex functional $c$. We assume that $c$ is presented by an *evaluation oracle* that given $x \in \mathbb{R}^d$ returns $c(x)$, and that $M$ is presented by an *independence oracle* that given $J \subseteq N$ asserts whether or not $J$ is an independent set of $M$. In this article we establish the following theorem.

**Theorem 1.4** *For any fixed $d$, the convex matroid optimization problem with oracle presented matroid $M$ over $N := \{1, \ldots, n\}$, weighting $w : N \longrightarrow \mathbb{R}^d$, and oracle presented convex functional $c : \mathbb{R}^d \longrightarrow \mathbb{R}$, can be solved in polynomial oracle time using $O(n^{2d-1} \log n)$ operations and queries.*

The computational complexity is measured in terms of the number of real arithmetic operations and oracle queries. For rational input the algorithm is (strongly) polynomial time in the Turing computation model where the input includes the binary encoding of the weighting $w : N \longrightarrow \mathbb{Q}^d$ and the binary encoding of an upper bound $U := \max_{J \subseteq N} c(w(J))$ on the relevant values of the convex functional, but we do not dwell on the details here.

The special case of the convex matroid optimization problem for *uniform matroids* coincides with the special case of the so-called *shaped partition problem* [7] for *two parts*. Therefore, the specializations to two-parts of the lower bounds of [1, 2] imply a lower bound of $\Omega(n^{d-1})$ on the complexity of the convex matroid optimization problem. It would be very interesting to further study a plausible common generalization of the convex matroid optimization problem for arbitrary matroids and the shaped partition problem for arbitrary number of parts.

## 2    Proof of the theorem

For a matroid $M = (N, \mathcal{B})$ and a weighting $w : N \longrightarrow \mathbb{R}^d$, consider the following convex polytope

$$\mathcal{P}_w^M \quad := \quad \text{conv} \{ w(B) \ : \ B \in \mathcal{B} \} \quad \subset \quad \mathbb{R}^d \ .$$

The convex matroid problem can be reduced to maximizing the convex functional $c$ over $\mathcal{P}_w^M$: there will always be an optimal basis $B \in \mathcal{B}$ for which $w(B)$ is a vertex of $\mathcal{P}_w^M$ and so the problem can be solved by picking the best such vertex. However, as the number of matroid bases is typically exponential in $n$ it is not possible to construct $\mathcal{P}_w^M$ directly in polynomial time. To overcome this we consider the following zonotope:

$$\mathcal{P}_w \quad := \quad \sum_{1 \leq i < j \leq n} [-1, 1] \cdot (w(i) - w(j)) \quad \subset \quad \mathbb{R}^d \quad .$$

**Proposition 2.1** *Fix any $d$. Then the number of vertices of the zonotope $\mathcal{P}_w$ is $O(n^{2(d-1)})$. Further, in polynomial time using that many arithmetic operations, all its vertices can be listed, each vertex $v$ along with a linear functional $a(v) \in \mathbb{R}^d$ uniquely maximized over $\mathcal{P}_w$ at $v$.*

*Proof.* The zonotope $\mathcal{P}_w$ is the Minkowski sum of $m := \binom{n}{2}$ line segments in $\mathbb{R}^d$ and therefore (cf. [5, 8]) has $O(m^{d-1}) = O(n^{2(d-1)})$ vertices which can all be enumerated, each $v$ along with a vector $a(v)$ uniquely maximized at $v$, using that many arithmetic operations. □

Let $\mathcal{P}^M := \mathrm{conv}\{\mathbf{1}_B : B \in \mathcal{B}\} \subset \mathbb{R}^n$ be the basis polytope of the matroid $M = (N, \mathcal{B})$, where $\mathbf{1}_B := \sum_{j \in B} e_j$ is the incidence vector of $B \in \mathcal{B}$ with $e_j$ the $j$th standard unit vector in $\mathbb{R}^n$. We include the short proof of the following statement.

**Proposition 2.2** *Every edge of the basis polytope is equal to $e_i - e_j$ for some pair $i, j \in N$.*

*Proof.* Consider any pair $A, B \in \mathcal{B}$ of bases such that $[\mathbf{1}_A, \mathbf{1}_B]$ is an edge (that is, a 1-face) of $\mathcal{P}^M$, and let $a \in \mathbb{R}^n$ be a linear functional uniquely maximized over $\mathcal{P}^M$ at that edge. If $A \setminus B = \{i\}$ is a singleton then $B \setminus A = \{j\}$ is a singleton as well in which case $\mathbf{1}_A - \mathbf{1}_B = e_i - e_j$ and we are done. Suppose then, indirectly, that it is not, and pick an element $i$ in the symmetric difference $A \triangle B := (A \setminus B) \cup (B \setminus A)$ of $A$ and $B$ of minimum value $a_i$. Without loss of generality assume $i \in A \setminus B$. Then there is a $j \in B \setminus A$ such that $C := A \setminus \{i\} \cup \{j\}$ is a basis of $M$. Since $|A \triangle B| > 2$, $C$ is neither $A$ nor $B$. By the choice of $i$, this basis satisfies $a \cdot \mathbf{1}_C = a \cdot \mathbf{1}_A - a_i + a_j \geq a \cdot \mathbf{1}_A$, and hence $\mathbf{1}_C$ is also a maximizer of $a$ over $\mathcal{P}^M$ so lies in the 1-face $[\mathbf{1}_A, \mathbf{1}_B]$. But no $\{0, 1\}$-vector is a convex combination of others, yielding a contradiction. □

The *normal cone* of a face at a polyhedron $P$ in $\mathbb{R}^d$ is the relatively open cone of those linear functionals $a \in \mathbb{R}^d$ uniquely maximized over $P$ at that face. The collection of normal cones of all faces of $P$ is called the *normal fan* of $P$. A polyhedron $P$ is a *refinement* of a polyhedron $Q$ if the normal fan of $P$ is a refinement of that of $Q$, that is, the closure of each normal cone of $Q$ is the union of closures of normal cones of $P$. We have the following lemma.

**Lemma 2.3** *The zonotope $\mathcal{P}_w$ is a refinement of the polytope $\mathcal{P}_w^M$.*

*Proof.* Let $\pi : \mathbb{R}^n \longrightarrow \mathbb{R}^d : e_j \mapsto w(j)$ be the natural projection sending the unit vector $e_j$ corresponding to the matroid element $j \in N$ to the vector $w(j) \in \mathbb{R}^d$. Then for each $B \in \mathcal{B}$ we have $\pi(\mathbf{1}_B) = w(B)$ and hence

$$\mathcal{P}_w^M \; = \; \mathrm{conv}\{w(B) : B \in \mathcal{B}\} \; = \; \mathrm{conv}\{\pi(\mathbf{1}_B) : B \in \mathcal{B}\} \; = \; \pi(\mathcal{P}^M)$$

so $\mathcal{P}_w^M$ is a projection of $\mathcal{P}^M$. Thus, each edge of $\mathcal{P}_w^M$ is the projection of some edge of $\mathcal{P}^M$ and hence, by Proposition 2.2, is equal to $\pi(e_i - e_j) = w(i) - w(j)$ for some pair $i, j \in N$. Thus, the zonotope $\mathcal{P}_w = \sum_{1 \le i < j \le n} [-1, 1] \cdot (w(i) - w(j))$ is the Minkowski sum of a set of segments containing all edge directions of $\mathcal{P}_w^M$ and hence its normal fan is a refinement of the normal fan of $\mathcal{P}_w^M$. □

We are now in position to prove our theorem.

*Proof of Theorem 1.4.* Given data $M, w, c$, the algorithm proceeds with the following steps: first, compute via Proposition 2.1 the list of $O(n^{2(d-1)})$ vertices $v$ of $\mathcal{P}_w$, each $v$ along with a linear functional $a(v) \in \mathbb{R}^d$ uniquely maximized over $\mathcal{P}_w$ at $v$. Second, for each $v$ do the following: let $a := a(v)$ and define the following weighting of matroid elements by scalars:

$$b : M \longrightarrow \mathbb{R} \ : \ j \mapsto a \cdot w(j) = \sum_{i=1}^d a_i w(j)_i \ ;$$

now apply a greedy algorithm to obtain a basis $B(v) \in \mathcal{B}$ of maximum weight $b(B)$, that is, sort $N$ by decreasing $b$-value (using $O(n \log n)$ operations) and find, using at most $n$ calls to the independence oracle presenting $M$, the lexicographically first basis $B(v)$. Third, for each $v$ compute the value $c(w(B(v)))$ using the evaluation oracle presenting $c$; an optimal basis for the convex matroid optimization problem is any $B(v)$ achieving maximal such value among the bases $B(v)$ of vertices $v$ of $\mathcal{P}_w$. The complexity is dominated by the second step which takes $O(n \log n)$ operations and queries and is repeated $O(n^{2(d-1)})$ times, giving the claimed bound.

We now justify the algorithm. First, we claim that each vertex $u$ of $\mathcal{P}_w^M$ satisfies $u = w(B(v))$ for some $B(v)$ produced in the second step of the algorithm. Consider any such vertex $u$. Since $\mathcal{P}_w$ refines $\mathcal{P}_w^M$ by Lemma 2.3, the normal cone of $u$ at $\mathcal{P}_w^M$ contains the normal cone of some (possibly more than one) vertex $v$ of $\mathcal{P}_w$. Then $a := a(v)$ is uniquely maximized over $\mathcal{P}_w^M$ at $u$. Now, consider the second step of the algorithm applied to $v$ and let $b$ be the corresponding scalar weighting of matroid elements. Then the $b$-weight of any basis $B$ satisfies

$$b(B) \ = \ \sum_{j \in B} a \cdot w(j) \ = \ a \cdot \sum_{j \in B} w(j) \ = \ a \cdot w(B) \ \le \ a \cdot u$$

with equality if and only if $w(B) = u$. Thus, the maximum $b$-weight basis $B(v)$ produced by the greedy algorithm will satisfy $u = w(B(v))$. Thus, as claimed, each vertex $u$ of $\mathcal{P}_w^M$ is obtained as $u = w(B(v))$ for some $B(v)$.

Now, since $c$ is convex, the maximum value $c(w(B))$ of any basis $B \in \mathcal{B}$ will occur at some vertex $u = w(B(v))$ of $\mathcal{P}_w^M = \text{conv}\{ w(B) \ : \ B \in \mathcal{B} \}$. Therefore, any basis $B(v)$ with maximum value $c(w(B(v)))$ is an optimal solution to the convex matroid optimization problem. The third step of the algorithm produces such a basis and so the algorithm is justified. □

Shmuel Onn
*Technion - Israel Institute of Technology, 32000 Haifa, Israel,*
   and
*University of California at Davis, Davis, CA 95616, USA.*
*email: onn@ie.technion.ac.il, onn@math.ucdavis.edu,*
   *http://ie.technion.ac.il/~onn*